\documentclass[preprint,3p,12pt,pdf]{elsarticle}

\usepackage{lineno,hyperref}
\modulolinenumbers[5]

\usepackage{epstopdf}

\usepackage{tikz}
\usetikzlibrary{arrows}

\usepackage{pst-node}
\usepackage{tikz-cd}

\usepackage[shellescape]{gmp}

\usepackage{mathrsfs}
\usepackage{amssymb}
\usepackage{amsfonts}
\usepackage{latexsym}
\usepackage{mathtools}
\usepackage{xcolor}
\usepackage{wrapfig}
\usepackage{floatflt}

\usepackage{mathtools}
\usepackage{extarrows}

\usepackage{graphicx}

\usepackage{subcaption}

\usepackage{endnotes}

\usepackage{hyperref}

\def\lb{\label}

\newcommand{\er}[1]{\textrm{(\ref{#1})}}

\newtheorem{theorem}{\bf Theorem}[section]

\def\a{\alpha}         
\def\b{\beta}          
\def\g{\gamma}         
\def\G{\Gamma}

\def\k{\kappa}

\def\ve{\varepsilon}           

    \def\R{{\mathbb R}}       
    \def\N{{\mathbb N}}   

\let\ge\geqslant                 

\def\iy{\infty}

\def\el2{\ell^{\,2}}             \def\1{1\!\!1}

\def\Im{\mathop{\mathrm{Im}}\nolimits}

\def\Re{\mathop{\mathrm{Re}}\nolimits}

\def\BBox{\hspace{1mm}\vrule height6pt width5.5pt depth0pt \hspace{6pt}}

\newtheorem{proposition}[theorem]{\bf Proposition}

\let\ge\geqslant

\newcommand{\ca}{\begin{cases}}
	\newcommand{\ac}{\end{cases}}
\newcommand{\ma}{\begin{pmatrix}}
	\newcommand{\am}{\end{pmatrix}}
\renewcommand{\[}{\begin{equation}}
	\renewcommand{\]}{\end{equation}}
\def\eq{\begin{equation}}
	\def\qe{\end{equation}}

\def\BBox{\hspace{1mm}\vrule height6pt width5.5pt depth0pt \hspace{6pt}}

\begin{document}

	\begin{frontmatter}

		\title{Complete left tail asymptotic for the density of branching processes in the Schr\"oder case}

		\date{\today}

		\author
		{Anton A. Kutsenko}
	
		\address{Mathematical Institute for Machine Learning and Data
Science, KU Eichst\"att--Ingolstadt, Germany; email: akucenko@gmail.com}
	
	\begin{abstract}
		For the density of Galton-Watson processes in the Schr\"oder case, we derive a complete left tail asymptotic series consisting of power terms multiplied by periodic factors.  
	\end{abstract}


	
\end{frontmatter}


{\section{Introduction}\lb{sec0}}

We consider a simple Galton-Watson branching process $Z_t$ in the supercritical case with the minimum family size $1$ - the so-called Schr\"oder case. The probability of the minimum family size is $0<p_1<1$. Note that the case of non-zero extinction probability can usually be reduced to the supercritical case with the help of Harris-Sevastyanov transformation, see \cite{B1}. We assume that the mean of offspring distribution $E<+\iy$. Then one may define the {\it martingal limit} $W=\lim_{t\to+\iy}E^{-t}Z_t$. It is proven in \cite{BB1} that $p(x)$, the density of $W$, has an asymptotic
\[\lb{000}
 p(x)=x^{\a}V(x)+o(x^{\a}),\ \ \ x\to+0,
\]
with explicit $\a=-1-\log_Ep_1$ and a continuous, positive, multiplicatively periodic function, $V$, with period $E$. Further references to this asymptotic are always based on the principal work \cite{BB1} and do not provide any formula for $V(x)$, see, e.g., the corresponding remark in \cite{FW} or \cite{S} devoted to the Schr\"oder case. Recently, some explicit expressions for $V(x)$ in terms of Fourier coefficients of $1$-periodic Karlin-McGregor function and some values of $\G$-function are given in \cite{K2}. The derivation is based on the complete asymptotic series for discrete relative limit densities of the number of descendants provided in \cite{K}. However, even for the first asymptotic term, the derivation is somewhat informal because the continuous {\it martingal limit} and the discrete distribution of the relative limit densities differ significantly. In particular, only the first term in both asymptotics has a common nature, all other terms are completely different. Now, I found beautiful and independent steps allowing me to obtain the complete expansion (not only the first term)
\[\lb{000a}
 p(x)=x^{\a}V_1(x)+x^{\a+\b}V_2(x)+x^{\a+2\b}V_3(x)+...,\ \ \ x>0,
\]
with the certain value $\a$ defined above, $\b=-\log_Ep_1>0$, and explicit multiplicatively periodic functions $V_1$, $V_2$, $V_3$, ..., see \er{012} and \er{013}. Finally, we obtain a very efficient representation of $V_n$ in terms of Fourier coefficients of Karlin-McGregor functions and some values of $\G$-function, see \er{206}. The explicit form of $V_1(x)$ presented in \cite{K2} significantly helped to determine this series. However, we provide an independent proof - the main result is formulated in Theorem \ref{T1}. Summarizing above, the asymptotic consists of periodic oscillations amplified by power-law multipliers. As mentioned in, e.g., \cite{DIL,DMZ,CG}, such type of behavior can be important in applications in physics and biology.

{\section{Main results}\lb{sec1}}

Let us recall some known facts about branching processes, including some recent results obtained in \cite{K}. The Galton--Watson process is defined by
\[\lb{001}
 X_{t+1}=\sum_{j=1}^{X_t}\xi_{j,t},\ \ \ X_0=1,\ \ \ t\in\N\cup\{0\},
\]
where all $\xi_{j,t}$ are independent and identically-distributed natural number-valued random variables with the probability-generating function
\[\lb{002}
 P(z):=\mathbb{E}z^{\xi}=p_0+p_1z+p_2z^2+p_3z^3+....
\]
For simplicity, we consider the case when $P$ is entire. In this case, the first moment
\[\lb{003}
 E=P'(1)=p_1+2p_2+3p_3+...<+\iy
\]
is automatically finite. A polynomial $P$ is common in practice, but the results discussed below can be extended to a wide class of non-entire $P$. As the Introduction mentions, we assume $p_0=0$, $p_1\ne0$. Another natural assumption is $p_1<1$, otherwise the case $p_1=1$ is trivial. Under these assumptions, we can define
\[\lb{004}
 \Phi(z)=\lim_{t\to\iy}p_1^{-t}\underbrace{P\circ...\circ P}_{t}(z),
\]
which is analytic at least for $|z|<1$. The function $\Phi$ satisfies the Scr\"oder-type functional equation
\[\lb{005}
 \Phi(P(z))=p_1\Phi(z),\ \ \ \Phi(0)=0,\ \ \ \Phi'(0)=1.
\] 
The function $\Phi$ has an inverse
\[\lb{006}
 \Phi^{-1}(z)=\k_1z+\k_2z^2+\k_3z^3+\k_4z^4+...,
\]
analytic in some neighborhood of $z=0$. The coefficients $\k_j$ can be determined by differentiating the corresponding Poincar\'e-type functional equation
\[\lb{007}
 P(\Phi^{-1}(z))=\Phi^{-1}(p_1z),\ \ \ \Phi^{-1}(0)=0,\ \ \ (\Phi^{-1})'(0)=1
\]
at $z=0$. In particular,
\[\lb{008}
 \k_1=1,\ \ \ \k_2=\frac{p_2}{p_1^2-p_1},\ \ \ \k_3=\frac{2p_2\k_2+p_3}{p_1^3-p_1},\ \ \ ....
\]
It is shown in, e.g., \cite{D1} or \cite{K2}, that the density $p(x)$, see \er{000}, can be computed by
\[\lb{009}
 p(x)=\frac1{2\pi}\int_{-\iy}^{+\iy}\Pi(\mathbf{i}y)e^{\mathbf{i}yx}dy=\frac1{2\pi\mathbf{i}}\int_{\g}\Pi(z)e^{zx}dz,
\]
where $\g$ is a modification of the original contour $\mathbf{i}\R$ discussed below and
\[\lb{010}
 \Pi(z):=\lim_{t\to+\iy}\underbrace{P\circ...\circ P}_{t}(1-\frac{z}{E^{t}}).
\]
This function is entire, satisfying another Poincar\'e-type functional equation
\[\lb{011}
 P(\Pi(z))=\Pi(Ez),\ \ \ \Pi(0)=1,\ \ \ \Pi'(0)=-1.
\]
Note that $E>1$, since $p_0=0$, $0<p_1<1$, and $\sum p_j=1$, see \er{003}. The existence of the limit in \er{010} follows from the standard facts that the linearization of $P$ at $1$ leads to the multiplication of the increment of the argument by $E$. In other words, one can use $P(1+w)=1+Ew+O(w^2)$ in \er{010} and see the convergence. The same reason applies to \er{004}. The details related to the theory of Schr\"oder and Poincar\'e-type functional equations are available in, e.g., \cite{M}.  Everything is ready to derive an expansion of $p(x)$ that, in turn, gives the complete left tail asymptotic series
\begin{align}\notag
 p(x) & =\frac1{2\pi\mathbf{i}}\int_{\g}\Phi^{-1}(\Phi(\Pi(z)))e^{zx}dz  =\sum_{n=1}^{+\iy}\frac{\k_n}{2\pi\mathbf{i}}\int_{\g}\Phi(\Pi(z))^ne^{zx}dz 
 \\
\ &=\sum_{n=1}^{+\iy}x^{{-1-n\log_Ep_1}}V_n(x),\lb{012}
\end{align}
see \er{009}, where
\[\lb{013}
 V_n(x)=\frac{\k_nx^{1+n\log_Ep_1}}{2\pi\mathbf{i}}\int_{\g}\Phi(\Pi(z))^ne^{zx}dz.
\]
Using \er{005}, \er{011}, and \er{013}, it is easy to check that all $V_n$ are multiplicatively periodic
\[\lb{014}
 V_n(\frac{x}{E})=\frac{\k_nx^{\frac{\ln Ep_1^n}{\ln E}}}{2\pi\mathbf{i} p_1^n}\int_{\g}\Phi(\Pi(z))^ne^{\frac{zx}E}\frac{dz}{E}=\frac{\k_nx^{\frac{\ln Ep_1^n}{\ln E}}}{2\pi\mathbf{i}}\int_{\frac{\g}{E}}\Phi(\Pi(z))^ne^{zx}dz=V_n(x).
\]
The integration contour $\g$ in \er{013} should be chosen so that $\Pi(\g)$ lies in the domain, where $\Phi$ is invertible, and $\g$ must be symmetric and connects $-\mathbf{i}\iy$ with $+\mathbf{i}\iy$. One of the choices is $\g=\mathbf{i}\R+\ve$ with sufficiently large $\ve>0$, since $\Pi(z)\to0$ for $\Re z\ge0$ and $z\to\iy$. For such type of contours, the second integral in \er{014} along $\g$ and along $\g/E$ is equal to the same value, due to the Cauchy integral theorem. The exact power-law decay of $\Pi(z)$, when $\Re z\ge0$ and $z\to\iy$, can be determined from the identity 
\[\lb{pi1}
 \Pi(z)=\Phi^{-1}(K(\log_E z)z^{\frac{\ln p_1}{\ln E}}),
\]
 where $K$ is $1$-periodic function, see \er{200}. Following \cite{D1}, integral \er{009} exists. We assume that it exists in a strong sense meaning $\Pi(iy)\to0$ for $y\to\pm\iy$ - this also means $\Phi(\Pi(iy))\to0$ for $y\to\pm\iy$. Hence, $K(z)$, see \er{200} is analytic in a neighborhood of $[x+\frac{{\bf i}\pi}{2\ln E},x+1+\frac{{\bf i}\pi}{2\ln E}]$ for some large $x\in\R$. Due to periodicity, $K(z)$ is analytic in some strip of a positive width, a neighborhood of the line $\R+\frac{{\bf i}\pi}{2\ln E}$. Taking into account the fact that $K(z)$ is also analytic in the strip $|\Im z|<\frac{\pi}{2\ln E}$, see \cite{K}, \cite{D1} and \cite{Du1}, we deduce that $K(z)$ is analytic and periodic in a larger strip $|\Im z|<k$ for some $k>\frac{\pi}{2\ln E}$. This means that the contour $\log_E\g$ lies strictly inside the domain of the definition of $K$, see Fig. \ref{fig0}, and, hence, $K(\log_E z)$ is bounded and smooth for $z\in\g$. Thus, $\Pi(z)=O(z^{\frac{\ln p_1}{\ln E}})$, see \er{pi1}, and if $\frac{\ln p_1}{\ln E}<-1$, integrals \er{013} converge absolutely. Recall that $p_1<1$ and $E>1$. Hence, $\frac{\ln p_1}{\ln E}<-1$ is not a rare case. In fact, it can be shown that \er{013} converges anyway because $\frac{\ln p_1}{\ln E}<0$ and decaying function $z^{\frac{\ln p_1}{\ln E}}$ is multiplied by the oscillatory factor $e^{zx}$ for $\Im z\to\pm\iy$, but it requires more cumbersome calculations.
 
\begin{figure}[h]
	\center{\includegraphics[width=0.7\linewidth]{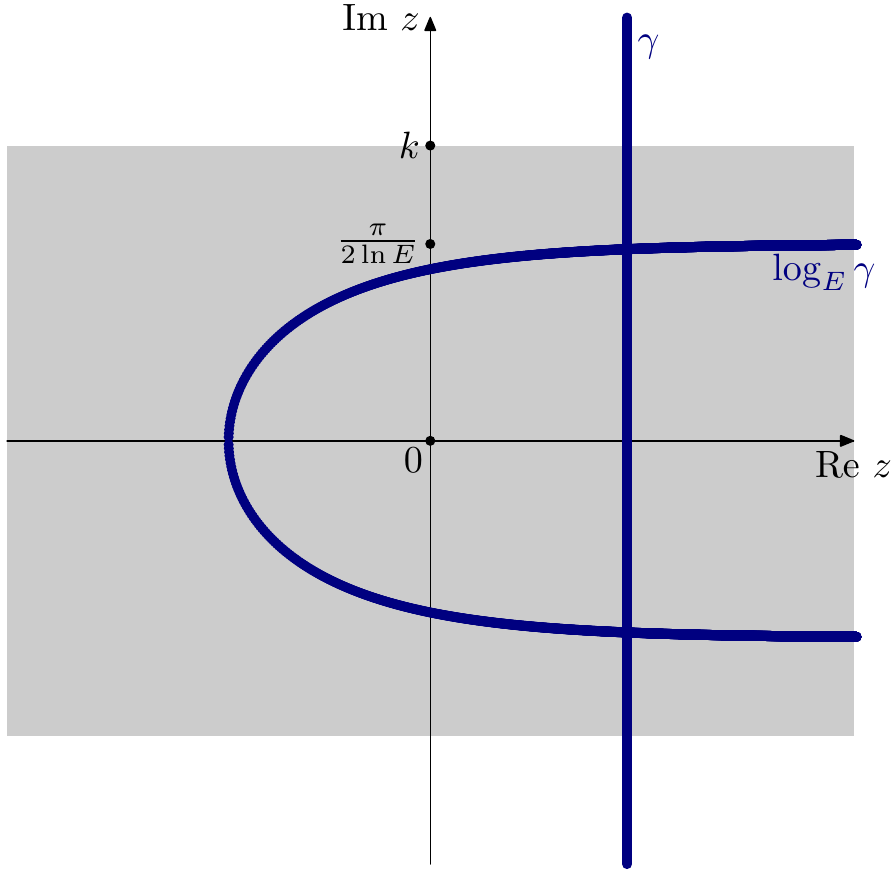}}
	\caption{Gray area is the domain of definition of $K(z)$. Paths $\gamma$ and $\log_E\gamma$ are plotted with the blue color.}\lb{fig0}
\end{figure}

\section{Representation of $V_n(x)$ through the Fourier coefficients of Karlin-McGregor function}

Recall that the Karlin-McGregor function, see \cite{KM1} and \cite{KM2}, is $1$-periodic function given by
\[\lb{200}
 K(z)=\Phi(\Pi(E^z))p_1^{-z}=\sum_{m=-\iy}^{m=+\iy}\theta_me^{2\pi\mathbf{i}mz},
\]
where the corresponding Fourier coefficients $\theta_m$ decay exponentially fast. The decay rate is determined by the width of the strip, where the Karlin-McGregor function is defined, see details in, e.g., \cite{K}. We recall only that the width is at least $\log_E\pi$. Let us rewrite \er{013} in terms of $K$:
\begin{multline}\lb{201}
 V_n(x)=\frac{\k_nx^{\frac{\ln Ep_1^n}{\ln E}}}{2\pi\mathbf{i}}\int_{\g}p_1^{-nz}\Phi(\Pi(E^z))^np_1^{nz}e^{E^zx}dE^z=\frac{\k_nx^{\frac{\ln Ep_1^n}{\ln E}}}{2\pi\mathbf{i}}\int_{\log_E\g}K(z)^np_1^{nz}e^{E^zx}dE^z\\
 =\frac{\k_nx^{\frac{\ln Ep_1^n}{\ln E}}}{2\pi\mathbf{i}}\int_{\log_E\g}(\sum_{m=-\iy}^{m=+\iy}\theta_{m}^{\ast n}e^{2\pi\mathbf{i}mz})p_1^{nz}e^{E^zx}dE^z=\k_n\sum_{m=-\iy}^{m=+\iy}\theta_{m}^{\ast n}\frac{x^{\frac{\ln Ep_1^n}{\ln E}}}{2\pi\mathbf{i}}\int_{\log_E\g}e^{2\pi\mathbf{i}mz}p_1^{nz}e^{E^zx}dE^z,
\end{multline} 
where $\theta_{m}^{\ast n}$ are Fourier coefficients of $1$-periodic function
\[\lb{202}
 K(z)^n=\sum_{m=-\iy}^{m=+\iy}\theta_{m}^{\ast n}e^{2\pi\mathbf{i}mz}.
\]
They are convolution powers of original Fourier coefficients $\theta_m$. We already know the form of $V_1(x)$:
\[\lb{203}
 V_1(x)=\sum_{m=-\iy}^{+\iy}\frac{\theta_me^{2\pi\mathbf{i}m\frac{-\ln x}{\ln E}}}{\Gamma(-\frac{2\pi\mathbf{i}m+\ln p_1}{\ln E})},
\]
see \cite{K2}. Comparing \er{201} with \er{203} and taking into account $\k_1=1$, see \er{008}, we guess the main identity
\[\lb{204}
 \frac{x^{\frac{\ln Ep_1}{\ln E}}}{2\pi\mathbf{i}}\int_{\log_E\g}e^{2\pi\mathbf{i}mz}p_1^{z}   e^{E^zx}dE^z=\frac{e^{2\pi\mathbf{i}m\frac{-\ln x}{\ln E}}}{\Gamma(-\frac{2\pi\mathbf{i}m+\ln p_1}{\ln E})}.
\]
A direct independent proof of \er{204} can be based on well-known Hankel integral representations of the $\G$-function. It is easy to check that, after changes of variables, \er{204} becomes equivalent to the formula presented in Example 12.2.6 on page 254 of \cite{WW}. This is a standard, but interesting, exercise in complex analysis and special functions. Substituting $p_1^n$ instead of $p_1$ into \er{204}, we obtain also
\[\lb{205}
 \frac{x^{\frac{\ln Ep_1^n}{\ln E}}}{2\pi\mathbf{i}}\int_{\log_E\g}e^{2\pi\mathbf{i}mz}p_1^{nz}e^{E^zx}dE^z=\frac{e^{2\pi\mathbf{i}m\frac{-\ln x}{\ln E}}}{\Gamma(-\frac{2\pi\mathbf{i}m+n\ln p_1}{\ln E})}.
\]
Combining \er{201} with \er{205}, we obtain the final result
\[\lb{206}
 V_n(x)=K_n(\frac{-\ln x}{\ln E}),\ \ \ K_n(z)=\k_n\sum_{m=-\iy}^{+\iy}\frac{\theta_{m}^{\ast n}e^{2\pi\mathbf{i}mz}}{\Gamma(-\frac{2\pi\mathbf{i}m+n\ln p_1}{\ln E})}.
\]
The formula \er{206} is more convenient for computations than \er{013} because it is similar to \er{203} for which efficient numerical schemes developed in \cite{K} were already applied in \cite{K2}. 

To justify changing the order of summation and integration in \er{201}, it is enough to take into account $\frac{\ln p_1}{\ln E}<-1$ and remember that $K$ is smooth and bounded in the symmetric strip, parallel to $\R$, of the width $\log_E\pi$, see details at the end of Section \ref{sec1}. Hence, $K(z)$ is well approximated by its Fourier series for $z\in\log_E\g$, with the necessary rate of convergence. Gathering all these remarks, let us formulate the main result.

\begin{theorem}\lb{T1} If $P$ is entire, integral \er{009} exists in the strong sense $\Pi(\mathbf{i}y)\to0$ for $y\to\pm\iy$, and $\frac{\ln p_1}{\ln E}<-1$ then \er{000a} along with \er{206} hold true.
\end{theorem}

As it is mentioned at the end of Section \ref{sec1}, the assumption $\frac{\ln p_1}{\ln E}<-1$ can be usually omitted. The convergence of \er{206} is quite fast because $\theta_{m}^{\ast n}=O(e^{-2\pi|m|k})$ for some $k>\frac{\pi}{2\ln E}$, which more than compensates small values of $\G$ in the denominator, see, e.g., the corresponding analysis based on Stirling approximations of Gamma function explained in \cite{K}. Thus, Fourier coefficients in \er{206} decay exponentially fast.

The condition $\Pi(\mathbf{i}y)\to0$ for $y\to\pm\iy$ can be replaced with a more simple condition, which follows directly from \er{011} and from the definition of the Julia set.
\begin{proposition}\lb{P1}
If $|\Pi(\mathbf{i}y)|<1$ for $y\in[r,Er]$, where $r>0$, then $\Pi(\mathbf{i}y)\to0$ for $y\to\pm\iy$.
\end{proposition}
{\it Proof.} Indeed, recall that the filled Julia set related to $P(z)$ contains the open unit disk because all the coefficients of the polynomial $P(z)$ are positive and the maximal value at the boundary of the unit disk is $P(1)=1$. Thus, if $\Pi(\mathbf{i}y)$ lies in the interior of the unit disk then $\Pi(\mathbf{i}y)$ lies in the interior of the corresponding component of the filled Julia set as well. Thus, $P$-iterations performed by the first formula in \er{011} converges to $0$ by the standard properties of the Julia sets, since $0$ is the unique attracting point inside the unit disc for $P(z)$. \BBox

Using fast exponentially convergent algorithms for the computation of $\Pi(z)$ developed in, e.g., \cite{K}, we can check the condition formulated in Proposition \ref{P1} numerically. On the other hand, if the interior of the corresponding component of the filled Julia set contains a little bit more than the open unit disk, namely a small sector near $z=1$ of the angle $\ge\pi$, then the condition formulated in Proposition \ref{P1} is automatically satisfied, since the point $w$, defined by
$$
w=\Pi(\mathbf{i}y)=1-\mathbf{i}y-\frac{\Pi''(0)}2y^2+O(y^3),\ \ \Pi''(0)=\frac{P''(0)}{E^2-E}>0,
$$ 
lies inside the filled Julia set and after a few, say $m$, $P$-iterations, see \er{011}, $\Pi(E^m\mathbf{i}y)=\underbrace{P\circ...\circ P}_{m}(w)$ will be small enough by the properties of Julia sets. All the filled Julia sets we tested in various examples contain such sectors, see, e.g., the right panel in Fig. \ref{fig2} computed with the help of \href{https://www.marksmath.org/visualization/polynomial_julia_sets/}{this site}\endnote{\it https://www.marksmath.org/visualization/polynomial\_julia\_sets/ \lb{ref1}}. 

\begin{figure}[h]
	\center{\includegraphics[width=0.99\linewidth]{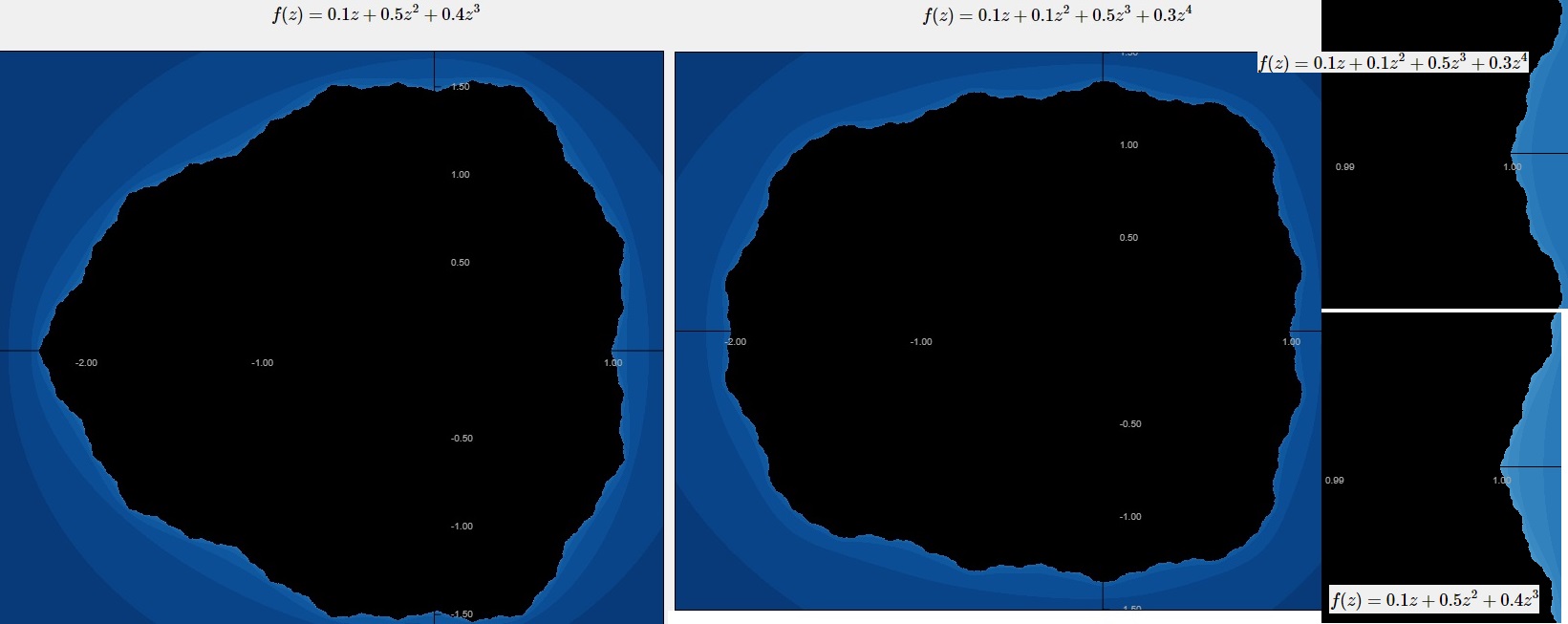}}
	\caption{Two examples of Julia sets  computed with the help of {\href{https://www.marksmath.org/visualization/polynomial_julia_sets/}{this site}$^1$}: the cases $p_1=0.1$, $p_2=0.5$, $p_3=0.4$, and $p_1=0.1$, $p_2=0.1$, $p_3=0.5$, $p_4=0.3$.}\lb{fig2}
\end{figure} 


\section{Examples}

As examples, we consider the same cases as in \cite{K2}, see Fig. \ref{fig1}. The difference is that we focus on the second asymptotic term. We use the same numerical procedures based on fast algorithms developed in \cite{K} plus FFT to compute $p(x)$. At the moment, FFT is a bottleneck because it is applied to the infinite integral, see the first formula in \er{009}, with the use of some tricks. The tricks are necessary because FFT is developed for the finite integrals - from $-\pi$ to $\pi$. This leads to insufficient accuracy for small $x$. Now, I am working on the improvement of the situation. However, this fact shows how important the asymptotic is, because the computation of $V_1$ and $V_2$ is sufficiently accurate and does not involve infinite integrals. 

\begin{figure}[h]
	\center{\includegraphics[width=0.9\linewidth]{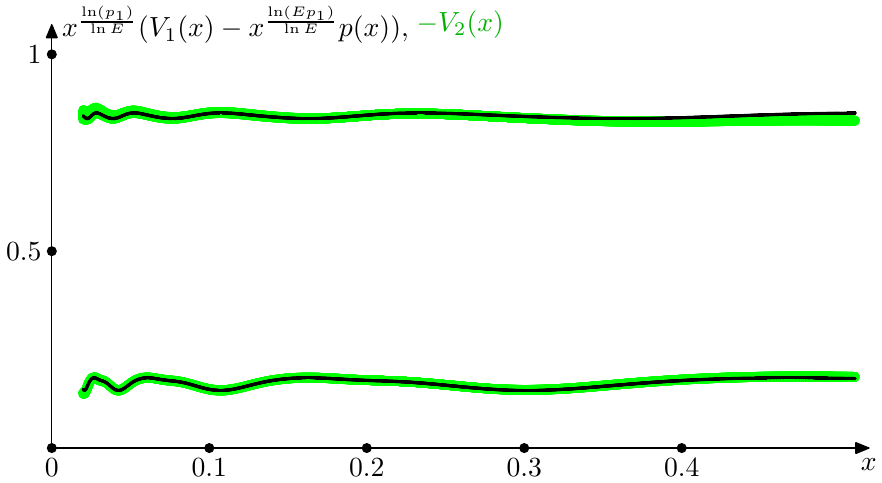}}
	\caption{Two examples: the cases $p_1=0.1$, $p_2=0.5$, $p_3=0.4$ (upper curves), and $p_1=0.1$, $p_2=0.1$, $p_3=0.5$, $p_4=0.3$.}\lb{fig1}
\end{figure} 

\section*{Acknowledgements} 
This paper is a contribution to the project M3 of the Collaborative Research Centre TRR 181 "Energy Transfer in Atmosphere and Ocean" funded by the Deutsche Forschungsgemeinschaft (DFG, German Research Foundation) - Projektnummer 274762653. 

%

\theendnotes

\end{document}